\documentclass[3p,times]{elsarticle}

\usepackage{amssymb}
\usepackage{amsmath}

\newcommand{\RR}{\mathbb{R}}

\newtheorem{theorem}{Theorem}
\newtheorem{corol}{Corollary}

\newtheorem{remark}{Remark}

\begin{document}

\begin{frontmatter}




\title{An efficient iterative technique for solving initial value problem for multidimensional partial differential equations}

\author{Josef Rebenda\fnref{auth1}}
\author{Zden\v{e}k \v{S}marda\corref{auth2}}
\fntext[auth1]{CEITEC BUT, Brno University of Technology, Technicka 3058/10, 61600 Brno, Czech Republic}
\cortext[auth2]{Department of Mathematics, Faculty of Electrical Engineering and Communication, Brno University of Technology, Technicka 8,  616 00 Brno, Czech Republic}

\address{}

\begin{abstract}
 
A new iterative technique is presented for solving of initial value problem for certain classes of multidimensional linear and nonlinear partial differential equations. Proposed iterative scheme does not require any discretization, linearization or small perturbations and therefore significantly reduces numerical computations. Rigorous convergence analysis of presented technique and an error estimate are included as well. Several numerical examples for high dimensional initial value problem for heat and wave type partial differential equations are presented to demonstrate reliability and performance of proposed iterative scheme.

\end{abstract}


\end{frontmatter}



\section{Introduction}
Many physical phenomena can be described by mathematical models that involve partial
differential equations. Therefore, in recent years, researchers look for new numerical methods which are more cost effective and
simple in implementation to solve partial differential equations. Investigation of exact and approximate
solution helps us to understand meaning and relevance of these mathematical models. Several techniques
including scattering method \cite{vak}, sine-cosine method \cite{wazwaz},  homotopy analysis method \cite{liao}, \cite{jafari}, homotopy perturbation method [5-9], differential transform method  \cite{solt}, variational iteration method \cite{he3}, \cite{biazar}, or decomposition methods [13-20] have been used for solving these problems, but mostly for two dimensional partial differential equations.

Inspired and motivated by ongoing research in this area, we apply new iterative scheme for solving heat- and
wave-type equations. Several examples are given to verify reliability and efficiency of proposed technique.

\section{Preliminaries}
Let $\Omega$ be a compact subset of $\RR^k$. Denote $J = [-\delta,\delta] \times \Omega$, where $\delta>0$ will be specified later, then $J$ is a compact subset of $\RR^{k+1}$. Let $u(t,x)=u(t,x_1,\dots,x_k)$
be a real function of $k+1$ variables defined on $J$. We introduce the following operators: $\nabla = (\frac{\partial}{\partial t},\frac{\partial}{\partial x_1}, \dots,
\frac{\partial}{\partial x_k})$ and $D = (\frac{\partial}{\partial x_1}, \dots, \frac{\partial}{\partial x_k})$. We deal with partial differential equations
\begin{equation}\label{eq1}
\frac{\partial^n}{\partial t^n} u(t,x) = F (t, x, u, \nabla u, \dots, \nabla^m u) \quad \text{for } m<n
\end{equation}
and
\begin{equation}\label{eq2}
\frac{\partial^n}{\partial t^n} u(t,x) = F (t, x, u, \nabla u, \dots, \nabla^{n-1} u, D \nabla^{n-1} u, D^2 \nabla^{n-1} u, \dots, D^{m-(n-1)} \nabla^{n-1} u) 
\quad \text{for } m \geq n.
\end{equation}
In both cases, left-hand side of the equation contains only the highest derivative with respect to $t$. We do not consider equations where the order of partial derivatives
with respect to $t$ is $n$ or higher on the righthand side, including mixed derivatives.

When convenient, we will use multiindex notation as well:
$$
\frac{\partial^{\vert \alpha \vert}}{\partial x^{\alpha}} = \frac{\partial^{\vert \alpha \vert}}{\partial x_1^{\alpha_1} \partial x_2^{\alpha_2} \dots \partial x_k^{\alpha_k}},
$$
where $\vert \alpha \vert = \alpha_1 + \alpha_2 + \dots + \alpha_k$.

Denote $N=\max \{ m,n \}$.
We consider equation \eqref{eq1} or \eqref{eq2} with the set of initial conditions
\begin{align}\label{eq3}
u(0,x) &= c_1 (x),\notag \\
\frac{\partial}{\partial t} u(0,x) &= c_2 (x), \notag \\
 &\vdots \notag \\
\frac{\partial^{n-1}}{\partial t^{n-1}} u(0,x) &= c_n (x),
\end{align}
where initial functions $c_i (x)$, $i=1,\dots,n$ are taken from space $C^N (\Omega, \RR)$. It means that we are looking for classical solutions.

For the purpose of clarity, we emphasize that our formulation covers for instance heat, wave, Burger, Boussinesq or Korteweg-de Vries (KdV) equations.

Obviously, $F \colon J \times \RR^k \times \RR \times \RR^{k+1} \times \RR^{(k+1)^2} \times \dots \times \RR^{(k+1)^m} \to \RR$ if $m<n$,
and $F \colon J \times \RR^k \times \RR \times \RR^{k+1} 
\times \RR^{(k+1)^2} \times \dots \times \RR^{(k+1)^{n-1}} \times \RR^{k(k+1)^{n-1}} \times \RR^{k^2(k+1)^{n-1}} \times \dots \times \RR^{k^{m-n+1}(k+1)^{n-1}} \to \RR$ if $m \geq n$.
Denote 
\begin{equation}
K_1 = \frac{(k+1)^{m+1} -1}{k} \quad \text{for } m<n
\end{equation}
and
\begin{equation}
K_2 = \frac{(k+1)^{n-1} -1}{k} + (k+1)^{n-1} \frac{k^{m-n+2} -1}{k-1} \quad \text{for } m \geq n.
\end{equation}
Then, if we consider $u$ as dependent variable, we see that $F$ is a function of $k+1+K_1$ variables in case $m<n$ or $k+1+K_2$ variables in case $m \geq n$.

Denote
\begin{equation}\label{eq4}
u_0 (t,x) = \sum_{i=1}^n c_i (x) \frac{t^{i-1}}{(i-1)!} = \sum_{i=1}^n \left( \frac{\partial^{i-1}}{\partial t^{i-1}} u(0,x) \right) \frac{t^{i-1}}{(i-1)!}.
\end{equation}
Then $u_0 \in C^N (J, \RR)$.

We suppose that $F$ is Lipschitz continuous in last $K_1$ ($m<n$) or $K_2$ ($m \geq n$) variables, i.e. $F$ satisfies condition
\begin{equation}\label{eq5}
\vert F (t,x,y_1, \dots, y_{K_l}) - F (t,x,z_1, \dots, z_{K_l}) \vert \leq L \left( \sum_{i=1}^{K_l} \vert y_i - z_i \vert \right), \quad l = 1 \text{ or } 2,
\end{equation}
on a compact set which is defined as follows: There is $R \in \RR, R>0$ such that \eqref{eq5} holds on 
\begin{equation}\label{eq6}
J \times \prod_{\alpha_0 + \vert \alpha \vert \leq m}
[c_{\alpha_0, \alpha}, d_{\alpha_0, \alpha} ],
\end{equation}
where
\begin{equation}\label{eq7}
c_{\alpha_0, \alpha} = \min_{(t,x) \in J} \left[ \frac{\partial^{\alpha_0 + \vert \alpha \vert}}{\partial t^{\alpha_0} \partial x^\alpha} u_0 (t,x) \right] -R,
\quad
d_{\alpha_0, \alpha} = \max_{(t,x) \in J} \left[ \frac{\partial^{\alpha_0 + \vert \alpha \vert}}{\partial t^{\alpha_0} \partial x^\alpha} u_0 (t,x) \right] +R,
\end{equation}
and $\alpha_0 <n$ in all cases.

Since $F$ is continuous on compact set, $\vert F \vert$ attains its maximal value on this set, denote it $M$. Then we put
\begin{equation}\label{eq10}
\delta = \left( \frac{R \cdot (n-1)!}{M} \right)^{1/n}.
\end{equation}


\section{Main results}

\begin{theorem}\label{th1}
Let the condition \eqref{eq5} hold. Then problem consisting of equation \eqref{eq1} or \eqref{eq2} and initial conditions \eqref{eq3} has a unique local solution on $(-\delta, \delta) \times \Omega$, where $\delta$ is defined by \eqref{eq10}.
\end{theorem}

{\bf Proof.} First define the following operator
\begin{equation}\label{eq8}
T u(t,x) = u_0 (t,x) + \int_0^t \frac{(t-\xi)^{n-1}}{(n-1)!} F \bigl(\xi, x, u (\xi,x), \nabla u(\xi,x), \dots \bigr) d \xi,
\end{equation}
where function $F$ has either $k+1+K_1$ or $k+1+K_2$ arguments and the last $K_1$ or $K_2$ arguments involve dependent variable $u$.

Starting with equation \eqref{eq1}, respective \eqref{eq2}, and using repeated integration by parts, it can be easily proved that if $u$ is a solution of equation $u=Tu$, i.e., if $u$ is a fixed point of operator $T$, then it is a solution of problem \eqref{eq1}, \eqref{eq3}, respective \eqref{eq2}, \eqref{eq3}.

Denote $J_1 = [-\delta_1, \delta_1] \times \Omega$, where $0<\delta_1<\delta$. Then $J_1$ is compact. Let $C^N (J_1, \RR)$ be the space of functions from $J_1$ to $\RR$ with continuous partial derivatives up to order $N$. This space is a Banach space with respect to the norm
\begin{equation}\label{eq9}
\Vert u \Vert_{C^N} = \sum_{\alpha_0 + \vert \alpha \vert \leq N} \max_{(t,x) \in J_1} \left| \frac{\partial^{\alpha_0 + \vert \alpha \vert}}{\partial t^{\alpha_0} \partial x^\alpha} u(t,x) \right|.
\end{equation}
It is obvious that, considering space $C^N$, the order of partial derivatives with respect to $t$ is allowed to be greater than or equal to $n$ in case $m \geq n$ when calculating the norm. Further, we define closed ball $B_R (u_0) \subseteq C^N (J_1,\RR)$ as follows:
$$
B_R (u_0) = \left\{ y \in C^N (J_1,\RR) \colon \Vert y-u_0 \Vert_{C^N} \leq R \right\}.
$$
It is not difficult to verify that $F$ composed with any $y \in B_R (u_0)$ and its appropriate derivatives satisfies Lipschitz condition \eqref{eq5} and that upper bound $\vert F \vert \leq M$ remains valid as well. Indeed, $J_1 \subseteq J$ and if $\Vert y-u_0 \Vert_{C^N} \leq R$, then
\begin{equation*}
\left\vert \frac{\partial^{\alpha_0 + \vert \alpha \vert}}{\partial t^{\alpha_0} \partial x^\alpha} \Bigl( y(t,x)-u_0 (t,x) \Bigr) \right\vert \leq R
\end{equation*}
for all $(t,x) \in J_1$ and for all $\alpha_0 + \vert \alpha \vert \leq N$. Hence point $\bigl(t, x, y(t,x), \nabla y(t,x), \dots, \nabla^m y(t,x) \bigr)$ or\\
$\bigl(t, x, y(t,x), \nabla y(t,x), \dots, D^{m-(n-1)} \nabla^{n-1} y(t,x) \bigr)$ respectively lies in the set $J \times \prod\limits_{\alpha_0 + \vert \alpha \vert \leq m}
[c_{\alpha_0, \alpha}, d_{\alpha_0, \alpha} ]$ defined by \eqref{eq6} and \eqref{eq7}.

We need to show that $T$ is a contraction on $B_R (u_0)$ for sufficiently small $\delta_1$.

In the first step, we prove that $T$ maps $B_R (u_0)$ into itself. Take any $y \in B_R (u_0)$. Then
\begin{align*}
\Vert Ty - u_0 \Vert_{C^N} &= \Bigl\Vert u_0 + \int_0^t \frac{(t-\xi)^{n-1}}{(n-1)!} F \bigl(\xi, x, y (\xi,x), \nabla y(\xi,x), \dots \bigr) d \xi - u_0 \Bigr\Vert_{C^N} \\
&\leq \Bigl\Vert \int_0^t \frac{t^{n-1}}{(n-1)!} F \bigl(\xi, x, y (\xi,x), \nabla y(\xi,x), \dots \bigr) d \xi \Bigr\Vert_{C^N} \leq \frac{\delta_1^{n-1}}{(n-1)!} M \Bigl\Vert \int_0^t d \xi \Bigr\Vert_{C^N} \\
&\leq \frac{\delta_1^n}{(n-1)!} M < R \ \ \text{for } \ 0<\delta_1 < \left( \frac{R \cdot (n-1)!}{M} \right)^{1/n}.
\end{align*}
Thus $Ty \in B_R (u_0)$ for $0<\delta_1 < \left( \frac{R \cdot (n-1)!}{M} \right)^{1/n}$.

The second step is to show that $T$ is a contraction. Choose arbitrary $y, z \in B_R (u_0)$. Then we have
\begin{align*}
\Vert Ty - Tz \Vert_{C^N} &= \Bigl\Vert \int_0^t \frac{(t-\xi)^{n-1}}{(n-1)!} \Bigl[ F \bigl(\xi, x, y (\xi,x), \nabla y(\xi,x), \dots \bigr) - F \bigl(\xi, x, z (\xi,x), \nabla z(\xi,x), \dots \bigr) \Bigr] d \xi \Bigr\Vert_{C^N} \\
&\leq \Bigl\Vert \frac{t^{n-1}}{(n-1)!} \int_0^t \bigl\vert F \bigl(\xi, x, y (\xi,x), \nabla y(\xi,x), \dots \bigr) - F \bigl(\xi, x, z (\xi,x), \nabla z(\xi,x), \dots \bigr) \bigr\vert d \xi \Bigr\Vert_{C^N}\\
&\leq \frac{\delta_1^{n-1}}{(n-1)!} \Bigl\Vert \int_0^t L \cdot \Biggl( \sum_{\alpha_0 + \vert \alpha \vert \leq m} \Biggl\vert \frac{\partial^{\alpha_0 + \vert \alpha \vert}}{\partial t^{\alpha_0} \partial x^\alpha} \Bigl( y(t,x)-z(t,x) \Bigr) \Biggr\vert \Biggr) d \xi \Bigr\Vert_{C^N}\\
&\leq \frac{\delta_1^{n-1}}{(n-1)!} L \ \Biggl\Vert \Biggl( \sum_{\alpha_0 + \vert \alpha \vert \leq m} \max_{(t,x) \in J_1} \Biggl\vert \frac{\partial^{\alpha_0 + \vert \alpha \vert}}{\partial t^{\alpha_0} \partial x^\alpha} \Bigl( y(t,x)-z(t,x) \Bigr) \Biggr\vert \Biggr) \int_0^t d \xi \Biggr\Vert_{C^N} \leq \frac{\delta_1^n}{(n-1)!} L \ \Vert y-z \Vert_{C^N},
\end{align*}
where $\alpha_0 <n$ for equation \eqref{eq2} and $L$ is a Lipschitz constant for $F$ introduced in \eqref{eq5}. It follows that $T$ is a contraction for $0<\delta_1 < \left( \frac{(n-1)!}{2L} \right)^{1/n}$.

Combining all results, we obtain that if $0<\delta_1 \leq \min \left\{ \frac{\delta}{2}, \left( \frac{R \cdot (n-1)!}{2M} \right)^{1/n}, \left( \frac{(n-1)!}{2L} \right)^{1/n} \right\}$, then operator $T$ is a contraction on $B_R (u_0)$. Applying Banach contraction principle, we can conclude that $T$ has a unique fixed point in $B_R (u_0)$ which is a unique solution of problem \eqref{eq1}, \eqref{eq3}, respective \eqref{eq2}, \eqref{eq3}.

Since $\delta_1$ depends only on the Lipschitz constant $L$ and on the distance $R$ from initial data to the boundaries of the intervals $[c_{\alpha_0, \alpha}, d_{\alpha_0, \alpha} ]$ wherein the estimate $M$ holds, we can 
apply our result repeatedly to get a unique local solution defined for $(t,x) \in (-\delta, \delta) \times \Omega$.
\qed

\begin{theorem}\label{th2}
Assume that condition \eqref{eq5} holds. Then iterative scheme $u_p = T u_{p-1}$, $p \geq 1$ with initial approximation $u_0$ defined by \eqref{eq4}, where $T$ is defined by \eqref{eq8}, converges to unique local solution $u(t,x)$ of problem \eqref{eq1}, \eqref{eq3}, respective \eqref{eq2}, \eqref{eq3}. Moreover, we have the following error estimate for this scheme:
\begin{equation}\label{eq13}
\Vert u(t,x) - u_p (t,x)\Vert_{C^N} \leq \frac{R \cdot \gamma^p}{1-\gamma}
\end{equation}
on $(-\delta_1, \delta_1) \times \Omega$, where $\delta_1$ is chosen such that operator $T$ is a contraction,
\begin{equation}\label{eq14}
\gamma = \frac{L \cdot \delta_1^n}{(n-1)!},
\end{equation}
and constants $L$ and $R$ are defined by \eqref{eq5} and \eqref{eq6}.
\end{theorem}

{\bf Proof.} First we need to show that sequence $\bigl( u_p \bigr)_{p=0}^\infty$ is convergent. We prove it by showing that it is a Cauchy sequence.
Take any $p,q \in \mathbb N$, $q \geq p$. Then
\begin{align*}
\Vert u_q - u_p \Vert_{C^N} &= \Bigl\Vert \int_0^t \frac{(t-\xi)^{n-1}}{(n-1)!}  \Bigl[ F \bigl(\xi, x, u_{q-1} (\xi,x), \nabla u_{q-1}(\xi,x), \dots \bigr) -  F \bigl(\xi, x, u_{p-1} (\xi,x), \nabla u_{p-1}(\xi,x), \dots \bigr) \Bigr] d \xi \Bigr\Vert_{C^N} \\
&\leq \frac{\delta_1^n}{(n-1)!} L \ \Vert u_{q-1}-u_{p-1} \Vert_{C^N} \leq \gamma \ \Vert u_{q-1}-u_{p-1} \Vert_{C^N},
\end{align*}
where $ 0 \leq \gamma \leq \frac{1}{2} <1$. Put $q=p+1$. We obtain
\begin{equation*}
\Vert u_{p+1} - u_p \Vert_{C^N} \leq \gamma \Vert u_{p}-u_{p-1} \Vert_{C^N} \leq \gamma^2 \Vert u_{p-1}-u_{p-2} \Vert_{C^N} \leq \ldots \leq \gamma^p \Vert u_1-u_0 \Vert_{C^N}.
\end{equation*}
Now, using the triangle inequality, we get
\begin{align*}
\Vert u_q - u_p \Vert_{C^N} &= \Vert u_q - u_{q-1} \Vert_{C^N} + \ldots + \Vert u_{p+2} - u_{p+1} \Vert_{C^N} + \Vert u_{p+1} - u_p \Vert_{C^N} \leq \bigl( \gamma^{q-1} +  \ldots + \gamma^{p+1} + \gamma^{p} \bigr) \Vert u_1-u_0 \Vert_{C^N} \\
&\leq \gamma^p \bigl( 1 + \gamma + \gamma^2 + \ldots + \gamma^{q-p-1} \bigr) \Vert u_1-u_0 \Vert_{C^N} \leq \gamma^p \frac{1-\gamma^{q-p}}{1-\gamma} \Vert u_1-u_0 \Vert_{C^N}.
\end{align*}
Since $\gamma <1$, then $1-\gamma^{q-p} <1$ as well, and we estimate
\begin{align}
\Vert u_q - u_p \Vert_{C^N} &\leq \frac{\gamma^p}{1-\gamma} \Vert u_1-u_0 \Vert_{C^N} \leq \frac{\gamma^p}{1-\gamma} \Vert Tu_0-u_0 \Vert_{C^N} = \frac{\gamma^p}{1-\gamma} \Bigl\Vert \int_0^t \frac{(t-\xi)^{n-1}}{(n-1)!} F \bigl(\xi, x, u_0 (\xi,x), \nabla u_0(\xi,x), \dots \bigr) d \xi \Bigr\Vert_{C^N} \notag \\
&\leq \frac{\gamma^p}{1-\gamma} \cdot \frac{\delta_1^{n}}{(n-1)!} \cdot M < \frac{\gamma^p}{1-\gamma} \cdot R \leq \frac{\left( \frac{1}{2} \right)^p}{\frac{1}{2}} \cdot R. \label{eq15}
\end{align}
It follows that for arbitrary $\epsilon>0$ there is a $P \in \mathbb N$, $P > 1- \log_2 \epsilon + \log_2 R$ such that if $p,q \geq P$, then $\Vert u_q - u_p \Vert_{C^N} < \epsilon$. Thus sequence $\bigl( u_p \bigr)_{p=0}^\infty$ is a Cauchy sequence and consequently a convergent sequence.

Then there is a limit $u=\lim\limits_{p \to \infty} u_p$ such that $u=Tu$. Hence $u$ is a fixed point of operator $T$. Applying Theorem \ref{th1} we conclude that this fixed point is unique and it is a unique local solution of problem \eqref{eq1}, \eqref{eq3}, respective \eqref{eq2}, \eqref{eq3}.

Error estimate \eqref{eq13} follows immediately from \eqref{eq15} by taking a limit for $q \to \infty$.
\qed

\begin{corol}\label{c1}
Let condition \eqref{eq5} be valid and suppose that $F$ can be written as $G+g$:
\begin{equation}\label{eq12}
F \bigl(t, x, z (t,x), \nabla z(t,x), \dots \bigr) = G \bigl(t, x, z (t,x), \nabla z(t,x), \dots \bigr) + g(t,x).
\end{equation}
Then we may choose initial approximation
\begin{equation}\label{eq11}
\bar{u}_0 = u_0 + \int_0^t \frac{(t-\xi)^{n-1}}{(n-1)!} g(\xi, x) d \xi = \sum_{i=1}^n \left( \frac{\partial^{i-1}}{\partial t^{i-1}} u(0,x) \right) \frac{t^{i-1}}{(i-1)!} + \int_0^t \frac{(t-\xi)^{n-1}}{(n-1)!} g(\xi, x) d \xi.
\end{equation}
\end{corol}

{\bf Proof.} Denote $M_2 = \max \{ \vert G \vert + \vert g \vert \}$. According to the proof of Theorem \ref{th1}, we only need to show that $\bar{u}_0 \in B_R (u_0)$, i.e., $\Vert \bar{u}_0 - u_0 \Vert_{C^N} \leq R$. Indeed, we have
\begin{align*}
\Vert \bar{u}_0 - u_0 \Vert_{C^N} &= \Bigl\Vert u_0 + \int_0^t \frac{(t-\xi)^{n-1}}{(n-1)!} g(\xi, x) d \xi  - u_0 \Bigr\Vert_{C^N} = \Bigl\Vert \int_0^t \frac{(t-\xi)^{n-1}}{(n-1)!} g(\xi, x) d \xi \Bigr\Vert_{C^N}\\
&\leq \Bigl\Vert \int_0^t \frac{(t-\xi)^{n-1}}{(n-1)!} \bigl\vert g(\xi, x) \bigr\vert d \xi \Bigr\Vert_{C^N} \leq \frac{t^{n-1}}{(n-1)!} \Bigl\Vert \int_0^t \bigl\vert g(\xi, x) \bigr\vert + \bigl\vert G \bigl(\xi, x, z (\xi,x), \nabla z(\xi,x), \dots \bigr)  \bigr\vert d \xi \Bigr\Vert_{C^N}\\
&\leq \frac{\delta_2^{n-1}}{(n-1)!} M_2  \Bigl\Vert \int_0^t d \xi \Bigr\Vert_{C^N} \leq \frac{\delta_2^n}{(n-1)!} M_2 < R \ \ \text{for } \ 0<\delta_2 < \left( \frac{R \cdot (n-1)!}{M_2} \right)^{1/n}.
\end{align*}
Consequently, integral
\begin{equation*}
\int_0^t \frac{(t-\xi)^{n-1}}{(n-1)!} g(\xi, x) d \xi
\end{equation*}
is small enough for sufficiently small $\delta_2$ and thus $\bar{u}_0 \in B_R (u_0)$ for this $\delta_2$.
\qed

\section{Applications}
We demonstrate potentiality of our approach on several initial value problems (IVP's).

{\bf Example 1.}  Consider the following two-dimensional heat-type equation 
\begin{equation}\label{p}
\frac{\partial u(x,y,t)}{\partial t} - \frac{\partial^2 u(x,y,t)}{\partial x^2}+ \frac{\partial^2 u(x,y,t)}{\partial y^2} + u(x,y,t) = (1+t) \sinh{(x+y)}
\end{equation}
with initial condition
\begin{equation}\label{p1}
u(x,y,0)= \sinh{(x+y)}.
\end{equation}
Then
$$
u_0(x,y,t) =\sum\limits_{k=0}^{n-1} {\frac{\partial ^ku\left({x,y,0} \right)}{\partial t^k}} \frac{t^k}{k!} + \int_0^t g(x,y,\xi) d\xi 
=  \sinh{(x+y)} \left(1 + t + \frac{t^2}{2}\right)
$$
and
$$
u_p(x,y,t) = u_0(x,y,t)+\int_0^t \left(\frac{\partial^2u_{p-1}(x,y,\xi)}{\partial x^2}- \frac{\partial^2 u(x,y,\xi)}{\partial y^2}- u(x,y,\xi)\right) d\xi, \ \ p \geq 1.
$$
Hence
\begin{eqnarray}
u_1(x,y,t) &=& \sinh{(x+y)} \left(1 + t + \frac{t^2}{2}\right)-\int_0^t \sinh{(x+y)} \left(1 + \xi + \frac{\xi^2}{2}\right) d\xi \nonumber \\
&=& \sinh (x+y)\left( 1- \frac{t^3}{3!}\right), \nonumber \\
u_2(x,y,t) &=& \sinh{(x+y)} \left(1 + t + \frac{t^2}{2}\right) - \int_0^t \sinh{(x+y)}\left( 1- \frac{\xi^3}{3!}\right)d\xi\nonumber \\
 &=& \sinh{(x+y)}\left(1+ \frac{t^2}{2}+ \frac{t^4}{4!}\right),\nonumber \\
u_3(x,y,t) &=&\sinh{(x+y)} \left(1 + t + \frac{t^2}{2}\right)- \int_0^t \sinh{(x+y)}\left(1+ \frac{\xi^2}{2}+  \frac{\xi^4}{4!}\right)d\xi \nonumber \\
&=& \sinh{(x+y)}\left( 1 + \frac{t^2}{2}-\frac{t^3}{3!}- \frac{t^5}{5!}\right),\nonumber \\  
u_4(x,y,t) &=& \sinh{(x+y)} \left(1 + t + \frac{t^2}{2}\right)- \int_0^t \sinh{(x+y)}\left( 1+ \frac{\xi^2}{2}-\frac{\xi^3}{3!}- \frac{\xi^5}{5!}\right)d\xi \nonumber \\
&=& \sinh{(x+y)}\left(1+ \frac{t^2}{2}-\frac{t^3}{3!} +\frac{t^4}{4!}+ \frac{t^6}{6!}\right),\nonumber \\  
u_5(x,y,t) &=&\sinh{(x+y)} \left(1 + t + \frac{t^2}{2}\right)- \int_0^t \sinh{(x+y)}\left( 1+ \frac{\xi^2}{2}-\frac{\xi^3}{3!}+\frac{\xi^4}{4}+ \frac{\xi^6}{6!} \right) d\xi \nonumber \\
&=& \sinh{(x+y)}\left( 1+\frac{t^2}{2}-\frac{t^3}{3!} +\frac{t^4}{4!}-\frac{t^5}{5!}-  \frac{t^7}{7!}\right),\nonumber \\  
\vdots \nonumber
\end{eqnarray}
We can see that so-called self-canceling terms appear between various components  (see, for example, $u_1,u_2,u_3$). Keeping the remaining non-canceled terms, we have
$$
u(x,y,t) = \sinh (x+y)\left( 1 +  \frac{t^2}{2!} -  \frac{t^3}{3!} + \frac{t^4}{4!} -  \frac{t^5}{5!} + \dots \right) = ( t+ e^{-t} )\sinh (x+y),
$$ 
which is unique exact solution of IVP \eqref{p}, \eqref{p1}.\\[3mm]
{\bf Example 2.} Consider the following initial value problem for two-dimensional heat-type equation with variable coefficients
\begin{equation}\label{3}
\frac{\partial u(x,y,t)}{\partial t} = \frac{y^2}{2} \frac{\partial^2 u(x,y,t)}{\partial x^2}+ \frac{x^2}{2} \frac{\partial^2 u(x,y,t)}{\partial y^2}
\end{equation}
with initial condition
\begin{equation}\label{4}
u(x,y,0)= y^2.
\end{equation}
Then
$$
u_0(x,y,t) = y^2
$$
and
$$
u_p(x,y,t) = y^2+\int_0^t \left(\frac{y^2}{2} \frac{\partial^2 u_{p-1}(x,y,\xi)}{\partial x^2}+ \frac{x^2}{2} \frac{\partial^2 u_{p-1}(x,y,\xi)}{\partial y^2}\right)d\xi,\ \ p \geq 1.
$$
Thus
\begin{eqnarray}
u_1(x,y,t) &=& y^2+\int_0^t x^2 d\xi = y^2+x^2t, \nonumber \\
u_2(x,y,t) &=& y^2+\int_0^t (y^2\xi + x^2 ) d\xi =   y^2\left( 1+\frac{t^2}{2!}\right) +x^2t,  \nonumber \\
u_3(x,y,t) &=& y^2+ \int_0^t \left( y^2\xi+ x^2 +x^2\frac{\xi^2}{2}\right) d\xi = y^2\left( 1+\frac{t^2}{2!}\right)+x^2 \left( t+  \frac{t^3}{3!} \right), \nonumber \\
u_4(x,y,t) &=& y^2\left(1+\frac{t^2}{2!} +\frac{t^4}{4!}\right ) + x^2 \left( t+  \frac{t^3}{3!} \right), \nonumber \\
u_5(x,y,t) &=& y^2\left(1+\frac{t^2}{2!} +\frac{t^4}{4!}\right ) + x^2 \left( t+  \frac{t^3}{3!}+ \frac{t^5}{5!} \right), \nonumber \\
\vdots \nonumber \\ 
u_{2l-1}(x,y,t) &=& y^2\left(1+\frac{t^2}{2!} +\frac{t^4}{4!} + \dots + \frac{t^{2l-2}}{(2l-2)!}\right )+  x^2  \left( t+  \frac{t^3}{3!}+ \frac{t^5}{5!}+\dots +\frac{t^{2l-1}}{(2l-1)!}\right), \nonumber \\
u_{2l}(x,y,t) &=&  y^2\left(1+\frac{t^2}{2!} +\frac{t^4}{4!} + \dots + \frac{t^{2l}}{(2l)!}\right )+  x^2  \left( t+  \frac{t^3}{3!}+ \frac{t^5}{5!}+\dots +\frac{t^{2l-1}}{(2l-1)!}\right), \nonumber \\
\vdots \nonumber 
\end{eqnarray}
Hence unique solution of IVP \eqref{3}, \eqref{4} has the form
\begin{eqnarray}
u(x,y,t) &=& x^2\left( t+ \frac{t^3}{3!} + \frac{t^5}{5!}+ \dots + \frac{t^{2l-1}}{(2l-1)!}+ \dots \right) + y^2\left( 1+ \frac{t^2}{2!} + 
\frac{t^4}{4!}+ \dots + \frac{t^{2l}}{(2l)!}+ \dots \right) \nonumber \\
&=& x^2\cosh{t} +y^2\sinh{t}. \nonumber
\end{eqnarray}
{\bf Example 3.} Consider nonlinear wave-type equation
\begin{equation}\label{pa}
\frac{\partial^2 u(x,y,t)}{\partial t^2}= 2x^2+2y^2 +\frac{15}{2} x \left( \frac{\partial^2 u(x,y,t)}{\partial x^2}\right )^2 + \frac{15}{2}y \left( \frac{\partial^2 u(x,y,t)}{\partial y^2}\right)^2 
\end{equation}
with initial conditions
\begin{equation}\label{pa1}
u(x,y,0) = \frac{\partial u(x,y,0)}{\partial t} = 0.
\end{equation}
Then
$$ 
u_0(x,y,t) = \int_0^t (t-\xi)(2x^2+2y^2)d\xi = t^2(x^2+y^2)
$$
and
$$
u_p(x,y,t) = t^2(x^2+y^2)+ \frac{15}{2}\int_0^t (t-\xi)\left(x \left( \frac{\partial^2 u_{p-1}(x,y,\xi)}{\partial x^2}\right )^2 + y \left( \frac{\partial^2 u_{p-1}(x,y,\xi)}{\partial y^2}\right)^2 \right) d\xi,
$$
$p\geq 1.$ From here we obtain the following iterations:
\begin{eqnarray}
u_1(x,y,t)&=& t^2(x^2+y^2)+ 30 \int_0^t (t-\xi)(x\xi^4 + y\xi^4) d\xi = t^2(x^2+y^2)+ t^6 (x+y), \nonumber \\
u_2(x,y,t)&=& t^2(x^2+y^2)+ 30 \int_0^t (t-\xi)(x\xi^4 + y\xi^4) d\xi = t^2(x^2+y^2)+ t^6 (x+y), \nonumber \\
\vdots \nonumber \\
u_l(x,y,t) &=& t^2(x^2+y^2)+ t^6 (x+y), \ l\geq 2.
\end{eqnarray}
Then we obtain required unique solution of \eqref{pa}, \eqref{pa1} as
$$
u(x,y,t)= t^2(x^2+y^2) + t^6(x+y).
$$
{\bf Example 4.} Consider three-dimensional wave-type equation with variable coefficients
\begin{equation}\label{h}
\frac{\partial^2 u}{\partial t^2}- \frac{1}{2} \left( x^2 \frac{\partial^2 u}{\partial x^2}+ y^2\frac{\partial^2 u}{\partial y^2}+z^2\frac{\partial^2 u}{\partial z^2}\right) = x^2+y^2+z^2
\end{equation}
with initial conditions
\begin{equation}\label{h1}
u(x,y,z,t) =0, \quad \frac{\partial u(x,y,z,t)}{\partial t} = x^2+y^2-z^2.
\end{equation}
Then
$$
u_0(x,y,z,t) = t(x^2+y^2-z^2) + \int_0^t (t-\xi)(x^2+y^2+z^2)d\xi = t(x^2+y^2-z^2)+\frac{t^2}{2!}(x^2+y^2+z^2)
$$
and
\begin{eqnarray*}
u_p(x,y,z,t) &=&  t(x^2+y^2-z^2)+\frac{t^2}{2!}(x^2+y^2+z^2)  \\
&+& \frac{1}{2}\int_0^t(t-\xi)\left( x^2 \frac{\partial^2 u_{p-1}}{\partial x^2}+ y^2\frac{\partial^2 u_{p-1}}{\partial y^2}+z^2\frac{\partial^2 u{p-1}}{\partial z^2}\right) d\xi, \ p\geq 1.
\end{eqnarray*}
From here we get
\begin{eqnarray*}
u_1(x,y,z,t)&=& t(x^2+y^2-z^2)+\frac{t^2}{2!}(x^2+y^2+z^2) \\
&+& \frac{1}{2}\int_0^t (t-\xi)\left[ (x^2+y^2)(2\xi+\xi^2) +z^2(-2\xi + \xi^2 )\right ] d\xi  \\
&=& t(x^2+y^2-z^2)+\frac{t^2}{2!}(x^2+y^2+z^2)+\frac{t^3}{3!}(x^2+y^2-z^2) +\frac{t^4}{4!}(x^2+y^2+z^2)  \\
&=& (x^2+y^2)\left(t+\frac{t^2}{2!} + \frac{t^3}{3!}+\frac{t^4}{4!} \right) + z^2 \left(-t+\frac{t^2}{2!}-\frac{t^3}{3!}+\frac{t^4}{4!}\right),\\
u_2(x,y,z,t)&=& (x^2+y^2)\left(t+\frac{t^2}{2!} + \frac{t^3}{3!}+\frac{t^4}{4!}+\frac{t^5}{5!} +\frac{t^6}{6!}\right) + z^2 \left(-t+\frac{t^2}{2!}-\frac{t^3}{3!}+\frac{t^4}{4!}
-\frac{t^5}{5!} +\frac{t^6}{6!}\right),\\
 \vdots  \\
u_k(x,y,z,t) &=& (x^2+y^2)\left(t+\frac{t^2}{2!} + \frac{t^3}{3!}+ \dots +\frac{t^{2k+2}}{(2k+2)!}\right) + z^2 \left(-t+\frac{t^2}{2!}-\frac{t^3}{3!}+ (-1)^k \frac{t^{2k+2}}{(2k+2)!}\right),\\
\vdots 
\end{eqnarray*}
Thus
\begin{eqnarray}
u(x,y,z,t) &=& \lim_{k\rightarrow \infty} \bigg [(x^2+y^2) \left( t+  \frac{t^2}{2!}+ \frac{t^3}{3!}+\dots + \frac{t^{2(k+1)}}{(2(k+1)!}\right)  \nonumber \\
 &+& z^2 \left (-t+\frac{t^2}{2!} - \frac{t^3}{3!}+ \dots  -\frac{t^{2k+1}}{(2k+1)!} \right)\bigg] \nonumber \\
 &=& (x^2+y^2)e^t + z^2e^{-t} -(x^2+y^2+z^2) \nonumber
 \end{eqnarray}
 which is unique exact solution of IVP \eqref{h},\eqref{h1}.

\begin{remark} Some above mentioned initial value problems have been solved using homotopy analysis method, homotopy perturbation method or Adomian decomposition method respectively 
(see \cite{jafari},\cite{jin},\cite{wazwaz1},\cite{adomian2}, \cite{din}). However, in contrast to iterative technique proposed in this paper, these methods require complicated calculations of multidimensional integrals, high order derivatives or Adomian's or He's polynomials as it can be seen in cited papers.
\end{remark}

\section{Conclusion}
\begin{itemize}
\item
We conclude that iterative algorithm presented in this paper is a powerful and efficient analytical technique suitable for numerical approximation of a solution of initial problem for wide class of partial differential equations of arbitrary order.
\item
There is no need for calculating multiple integrals or derivatives, only one integration in each step is performed. Less computational work is demanded compared to other methods (Adomian decomposition method, variational iteration method, homotopy perturbation method, homotopy analysis method).
\item
Expected solution is a limit of a sequence of functions, in contrast to other frequently used methods where a sum of a functional series is considered. Consequently, the form of a solution can be immediately controlled in each step.
\item
All notations are carefully described and proofs are treated rigorously, compared to many recently presented algorithms and methods.
\item
Region and rate of convergence depend on Lipschitz constant for righthand side $F$.
\item
Using presented approach, we are able not only to obtain approximate solution, but even there is a possibility to identify unique solution of initial problem in closed form.
\item
A specific advantage of this technique over any purely numerical method is that it offers a smooth, functional form of the solution in each step.
\item
Another advantage is that using our approach we avoided discretization, linearization or perturbation of the problem.
\item
There is a possibility to reduce computational effort by combining presented algorithm with Laplace transform since there is a convolution integral inside the iterative formula.
\item
Finally, a subject of further investigation is to develop the presented technique for systems of PDE's, to find modifications for solving equations with deviating arguments and for other types of problems (e.g. BVP's).
\end{itemize}

\section{Acknowledgements}
\label{sec6}
The first author was supported by the project CZ.1.07/2.3.00/30.0039 of Brno
University of Technology. The work of the second author was realised in CEITEC - Central European Institute of Technology with research infrastructure supported by the project CZ.1.05/1.1.00/02.0068 financed from European Regional Development Fund and by the project FEKT-S-11-2-921 of Faculty of Electrical Engineering and Communication, Brno University of Technology.
This support is gratefully acknowledged.





\bibliographystyle{elsarticle-num}
\bibliography{<your-bib-database>}



\end{document}